\documentclass[11pt]{article}
\usepackage{amsmath,amsthm,amsfonts,amssymb,amscd, amsxtra}

\newcommand{\banacha}{X}

\newtheorem{theorem}{Theorem}
\newtheorem{lemma}[theorem]{Lemma}
\newtheorem{definition}{Definition}
\newtheorem{corollary}[theorem]{Corollary}
\newtheorem{proposition}[theorem]{Proposition}

\newtheorem{remark}{Remark}

\begin{document}
\title{\textbf{Local convergence analysis of  Inexact Newton  method with relative residual error tolerance  under majorant condition in Riemannian Manifolds}}

\author{  T. Bittencourt \thanks{IME/UFG,  CP-131, CEP 74001-970 - Goi\^ania, GO, Brazil (Email: {\tt
      tiberio.b@gmail.com}). This author was supported  by CAPES.}
\and 
O. P. Ferreira\thanks{IME/UFG,  CP-131, CEP 74001-970 - Goi\^ania, GO, Brazil (Email: {\tt
      orizon@ufg.br}). This author was supported  by
     CNPq Grants 302024/2008-5, 480101/2008-6 and 473756/2009-9,  PRONEX--Optimization(FAPERJ/CNPq) and FUNAPE/UFG.}
}
\date{September 18, 2013}

\maketitle

\begin{abstract}
{A local convergence analysis of   Inexact  Newton's method with relative residual error tolerance  for finding a singularity  of a differentiable vector field defined on a complete Riemannian manifold, based on  majorant principle, is presented in this paper.  We prove that under local assumptions, the inexact Newton method with a fixed relative residual error tolerance converges Q -linearly to a singularity of the  vector field under consideration. Using this result we show that the  inexact Newton method  to find a zero of an analytic vector field can be implemented with a fixed relative residual error tolerance.  In the absence of errors, our analysis retrieve the classical local theorem on the Newton method in Riemannian context.
\\} 

\noindent
{{\bf Keywords:}  Inexact  Newton's method, majorant principle, local convergence analysis, Riemannian manifold.}
\end{abstract}

\maketitle
\section{Introduction}\label{sec:int}
Newton's method and its variations, including the inexact Newton methods, are the most efficient methods known for solving nonlinear equations  in Banach spaces.   Besides its practical applications, Newton's method is  also a powerful theoretical tool with a wide range of applications in pure and applied mathematics, see \cite{Blum1998, DR2009, Gondzio2013,  Nesterov1994, Potra2005, S1994, Wang1999, Wy1996}. In particular,  Newton's method has been instrumental in the modern complexity analysis of the solution of polynomial or analytical   equations \cite{Blum1998, Smale1986}, linear and quadratic programming problems and linear semi-definite programming problems \cite{Gondzio2012,  Gondzio2013,   Nesterov1994, Potra2005}.   In all these applications, homotopy methods are combined with Newton's method, which helps the algorithm to keep track of the solution of a parametrized  perturbed version of the original problem. 

In classic Newton's method, a linear equation system is solved in each iteration which can be expensive and unnecessary when the problem size is large. Inexact Newton's method comes up to overcome such drawback and can effectively cut down the computational cost by solving the linear equations approximately, see \cite{Dembo1982, DR2013, Morini1999}.  It would be most desirable to have an \emph{a priori} prescribed residual error tolerance in the iterative solutions of linear system for computing the Inexact Newton steps, in order to avoid under-solving or over-solving the linear system in question. The advantage of working with an error tolerance on the residual rests in the fact that the exact Newton step need not to be know for evaluating this error, which makes this criterion attractive for practical applications, see  \cite{Gondzio2012,  Gondzio2013}.

Newton's method has been  extended to Riemannian manifolds with many different  purposes. In particular,   in the last few years,  a couple of papers have dealt with the issue of  convergence analysis of Newton's method  for finding a singularity  of a differentiable vector field defined on a complete Riemannian manifold, see   \cite{Alvarez2008, D2003,  FS2012, FerreiraSvaiter2002, Li2006,  LW2008, Li2009, S1994, W2011, Wang2009,  Wang2006, Wang2007, WL2012}. Extensions to Riemannian manifolds of  analyses of Newton's method under the $\gamma$-condition  was given in \cite{D2003, Li2006, LW2008, Li2009}. Although the local convergence analysis of Inexact Newton's method in Banach space with relative errors tolerance  in the residue~\cite{Chen2006, Dembo1982, Morini1999}  are well understood,   as far as we know,  the convergence analysis of the method in Riemannian  manifolds context under general local assumptions, assuming \emph{only} bounded relative residual errors,  is a new contribution of this paper.    It is worth to point out  that, for null error tolerance,  the analysis   presented merge in the usual local convergence analysis on Newton's method in Riemannian manifold under a majorant condition,  see~\cite{FS2012}.  In our analysis, the classical Lipschitz condition is relaxed  using a  majorant function which provides a clear relationship between the majorant function  and the vector field under consideration.   Moreover,  several unrelated previous results pertaining to Newton's method are unified (see \cite{D2003, Li2006, LW2008}), now  in the Riemannian context.

The organization of the paper is as follows. In Section \ref{sec:int}, the notations  and basic results used in the paper are presented. In Section \ref{sec:lkant}  the main result is stated and in Section \ref{sec:pr} some properties of the majorant function are established and the main relationships between the majorant function and the vector field used in the paper are presented.  In Section 5  the main result is proved and two applications of this result are given in Section 6. Some final remarks are made in Section 7.

\section{Notation and auxiliary results} \label{sec:int}
In this section we recall some notations, definitions and basic properties of Riemannian manifolds used throughout the paper,  they can be found, for example in \cite{DoCa92}  and \cite{La95}.

Throughout the paper, $\cal{M}$ is a smooth manifold and  $C^{1}(\cal M)$ is the class of all  continuously differentiable functions  on  $\cal M$. The space of vector fields on $\cal M $ is denoted by ${\mathcal X}(\cal M)$,   by $T_{p}{\cal M}$ we denote the tangent space of $\cal M$ at $p$ and by $T{\cal M}={\bigcup_{x\in {\cal M}}}\,T_{x}{\cal M}$ the {\it tangent bundle \/} of $\cal M$ .
Let $\cal M$ be endowed with a Riemannian metric $ {\langle} \cdot , \cdot {\rangle}$, with corresponding norm denoted by $\|\cdot\|$, so that $\cal M$ is now a {\it Riemannian manifold}. Let us recall that the metric can be used to define the length of a piecewise $C^{1}$ curve $\zeta :\, [a,b] \rightarrow {\cal M}$ joining $p$ to $q$, i.e., such that $\zeta(a)=p$ and $\zeta(b)=q$, by $l(\zeta) = \int_a^b \|\zeta^{\prime}(t)\| dt$. Minimizing this length functional over the set of all such curves we obtain a distance $d(p,q)$, which induces the original topology on $ M$. The open and closed  balls of radius $r>0$ centered at $p$ are  defined, respectively, as
$$
B_{r}(p):=\left\{q\in M:d(p,q) <r\right\}, \qquad {\overline{B}}_{r}(p):=\left\{q\in M:d(p,q) \leq r\right\}.
$$
Also the metric induces a map $f\in C^{1}(\cal M) \mapsto$ $\mbox{grad} f \in {\cal X}(\cal M)$, which associates to each $f$ its  {\it gradient} via the rule $\langle \mbox{grad}f,X \rangle =df(X)$, for all $X\in  {\mathcal X}(\cal M)$. The chain rule generalizes to this setting in the usual way: $(f\circ \zeta)^{\prime}(t)=\langle \mbox{grad}f(\zeta(t)), \zeta^{\prime}(t)\rangle $, for all  curves $\zeta\in C^{1}$.  Let $\zeta$ be a curve joining the points $p$ and $q$ in $\cal M$ and let
$\nabla$ be a Levi-Civita connection associated to $(\cal M,{\langle}, {\rangle})$. For each $t \in [a,b]$, $\nabla$ induces an isometry, relative to $ \langle \cdot , \cdot \rangle  $,
\begin{equation} \label{D:dtp}
\begin{split}
P_{\zeta,a,t} \colon T _{\zeta(a)} {\cal M} &\longrightarrow T _
{\zeta(t)} {\cal M}\\
v &\longmapsto P_{\zeta,a,t}\, v = V(t),
\end{split}
\end{equation}
where $V$ is the unique vector field on $\zeta$ such that
$ \nabla_{\zeta'(t)}V(t) = 0$ and $V(a)=v$,
the so-called {\it parallel translation} along $\zeta$ from $\zeta(a)$ to $\zeta(t)$. Note also that
\[
P_{\zeta,b_1,b_2}\circ P_{\zeta,a,b_1} = P_{\zeta,a,b_2}, \qquad P_{\zeta,b,a}  = {P_{\zeta,a,b} } ^{-1}.
\]
A vector field $V$ along $\zeta$ is said to be {\it parallel} if $\nabla_{\zeta^{\prime}} V=0$. If $\zeta^{\prime}$ itself is parallel, then we say that $\zeta$ is a {\it geodesic}. The geodesic equation $\nabla_{\ \zeta^{\prime}} \zeta^{\prime}=0$ is a second order nonlinear ordinary differential equation, so the geodesic $\zeta$ is determined by its position $p$ and velocity $v$ at $p$. It is easy to check that $\|\zeta ^{\prime}\|$ is constant. We say that $ \zeta $ is {\it normalized} if $\| \zeta^{\prime}\|=1$.  A geodesic $\zeta: [a, b]\to {\cal M}$ is said to be {\it minimal} if its length is equal the distance of its end points, i.e. $l(\zeta)=d(\zeta(a), \zeta(b))$.

A Riemannian manifold is {\it complete} if its geodesics are defined for any values of $t$. The Hopf-Rinow's
theorem asserts that if this is the case then any pair of points, say $p$ and $q$, in $\mathcal{M}$ can be joined by
a (not necessarily unique) minimal geodesic segment. Moreover, $({\cal M}, d)$ is a complete metric space
and bounded and closed subsets are compact.

The {\it exponential map} at $p,$ $\mbox{exp}_{p}:T_{p}  {\cal M} \rightarrow {\cal M} $ is defined by $\mbox{exp}_{p}v\,=\, \zeta _{v}(1)$, where $\zeta _{v}$ is the geodesic defined by its position $p$ and velocity $v$ at $p$ and  $ \,\zeta _{v}(t)\,=\,\mbox{exp}_{p}tv$ for any value of $t$.
For $p\in {\cal M}$, let
$$
r_{p}:=\sup\left\{ r >0 : {\mbox{exp}_{p}}_{|_{B_{r}(o_{p})}} \,\;\mbox{is a diffeomorphism} \right\},
$$
where $o_{p}$ denotes the origin of $T_{p}{\cal M}$ and $B_{r}(o_{p}):=\{v\in T_{p}{\cal {\cal M}} : \|v-o_{p}\|<r\}$.  Note that if $0<\delta < r_{p}$ then $\mbox{exp}_{p} B_{\delta}(o_{p})=B_{\delta}(p)$. The number  $r_{p}$ is called the {\it injectivity radius} of $\cal M$ at $p$.
\begin{definition} \label{def:kp}
Let $p\in {\cal M}$ and $r_p$ the radius of injectivity at $p.$ Define the quantity
$$
K_{p}:=\sup\left\{\frac{d(\exp_{q}u,\exp_{q}v)}{\|u-v\|}~:~  q\in B_{r_p}(p), ~ u,v\in T_{q}{\cal M}, ~ u\neq v, ~ \|v\|\leq r_p, ~ \|u-v\|\leq r_p\right\}.
$$
\end{definition}

\begin{remark}
The quantity $K_{p}$ measures how fast the geodesics spread apart in $\cal M.$ In particular, when  $u=0$ or more generally when $u$ and $v$ are on the same line through $o_{q},$
$$d(\exp_{q}u,\exp_{q}v)= \|u-v\|.$$
So  $K_{p}\geq 1$ for all $p\in {\cal M}$.  When $\cal M$ has non-negative sectional curvature, the geodesics spread apart less than the rays (\cite{DoCa92}, Chap. 5) so that
 $$d(\exp_{q}u,\exp_{q}v)\leq \|u-v\|.$$
As a consequence  $K_{p}=1$ for all $p\in {\cal M}.$  Finally it is worth mentioning that radii less than $r_p$ could be used as well (although this would require added notation such as $K_p(\rho)$
for  $r_p$). In this case,  the measure by which geodesics spread apart might decrease, thereby providing slightly stronger results so long as the radius was not too much less than $r_p$.
\end{remark}
Let $X$ be a $C^1$ vector field on $\cal M$. The  covariant derivative of $X$ determined by
the Levi-Civita connection $\nabla$ defines at each $p\in {\cal M}$ a linear map $\nabla X (p):T_p {\cal M} \to T_p {\cal M}$ given by
\begin{equation}\label{D:ddc}
 \nabla X(p) v:=\nabla_Y X (p),
\end{equation}
where $Y$ is a vector field such that $Y(p)=v$.
\begin{definition} \label{d:scd} Let $\cal M$ be a complete Riemannian manifold and  $Y_1, \ldots, Y_n$ be vector fields on $\cal M$. Then,  the $n$-th covariant derivative of $X$ with respect to $Y_1, \ldots, Y_n$ is defined inductively by
\[
\nabla^{2}_{{\{Y_1, Y_2\}}}X:=\nabla_{{Y_2}}\nabla_{Y_1}X, \qquad
 \nabla^{n}_{{\{Y_i\}_{i=1}^{n}}}X:=\nabla_{{Y_n}}( \nabla_{Y_{n-1}} \cdots \nabla_{Y_{1}}X).
\]
\end{definition}
\begin{definition} \label{d:scdp} Let $\cal M$ be a complete Riemannian manifold, and $p\in {\cal M}$. Then,  the $n$-th covariant derivative of $X$  at $p$ is the $n$-th multilinear map $\nabla^n X (p):T_p {\cal M}\times \ldots \times T_p {\cal M} \to T_p {\cal M}$ defined by
\[
\nabla^{n}X(p)(v_1, \dots, v_n):=\nabla^{n}_{{\{Y_i\}_{i=1}^{n}}}X(p),
\]
where $Y_1, \ldots, Y_n$ are vector fields on $\cal M$ such that $Y_1(p)=v_1, \ldots, Y_n(p)=v_n$.
\end{definition}
We remark that Definition~\ref{d:scdp} only depends on the $n$-tuple of vectors $(v_1, \ldots, v_n)$ since the covariant derivative is tensorial in each vector field $Y_i$.

\begin{definition} \label{d:norm}
 Let $\cal M$ be a complete Riemannian manifold and $p\in {\cal M}$. The norm of an
$n$-th multilinear map $A:T_p {\cal M}\times \ldots \times T_p {\cal M} \to T_p {\cal M}$  is defined by
  \[
 \|A\|=\sup \left\{ \|A(v_1, \dots, v_n) \| \;:\;\; v_1, \dots, v_n\in T_p {\cal M}, \,\|v_i\|=1, \, i=1, \ldots, n \right\}.
 \]
 In particular the norm of the $n$-th covariant derivative of $X$  at $p$ is given  by
 \[
 \|\nabla^{n}X(p)\|=\sup\left\{\|\nabla^{n}X(p)(v_1, \dots, v_n) \| \;:\;\; v_1, \dots, v_n\in T_p {\cal M}, \,\|v_i\|=1, \, i=1, \ldots, n \right\}.
 \]
\end{definition}

\begin{lemma}\label{le:TFC}
Let $\Omega$ be an open subset of $\cal M$, $X$ a
$C^1$
vector field defined on $\Omega$ and  $\zeta:[a, b]\to {\Omega}$  a $C^\infty$ curve. Then
\[
P_{\zeta,t,a} X(\zeta(t))= X(\zeta(a)) + \int_a^t
   P_{\zeta,s,a} \nabla X (\zeta(s)) \,\zeta'(s) \,ds, \qquad t\in [a, b].
\]
\end{lemma}
\begin{proof}
See \cite{FerreiraSvaiter2002}.
\end{proof}
\begin{lemma}\label{le:TFC2}
Let $\Omega$ be an open subset of $\cal M$, $X$ a
$C^2$
vector field defined on $\Omega$ and  $\zeta:[a, b]\to {\Omega}$  a $C^\infty$ curve. Then for all $Y\in {\mathcal X}(\cal M)$ we have that
\[
P_{\zeta,t,a} \nabla X(\zeta(t))\,Y(\zeta(t))= \nabla X(\zeta(a))Y(\zeta(a)) + \int_a^t
   P_{\zeta,s,a} \nabla^{2} X (\zeta(s)) (Y(\zeta(s)), \zeta'(s)) \,ds,  \qquad t\in [a, b].
\]
\end{lemma}
\begin{proof}
See \cite{Li2006}.
\end{proof}
\begin{lemma}[Banach's Lemma] \label{lem:ban}
Let $B$ be a linear operator  and  let $I_p$ be the identity operator in  $T _p M$. If  $\|B-I_p\|<1$  then $B$ is invertible and
$
\|B^{-1}\|\leq 1/\left(1-  \|B-I_p\|\right).
$
\end{lemma}
\begin{proof}
Under the hypothesis, it is easily shown that $B^{-1}=\sum_{i=0}^{\infty}(B-I_p)^{i}$ and hence $\|B^{-1}\|\leqslant \sum_{i=0}^{\infty}\|B-I_p\|^{i}=1/(1-\|(B-I_p)\|).$
\end{proof}
\section{Local analysis  for  Inexact Newton method } \label{sec:lkant}
Our goal is to prove in Riemannian manifold context  the following version  of  Inexact Newton  method with relative residual error tolerance  under majorant condition.
\begin{theorem}\label{th:nt}
Let $\cal M$ be a Riemannian manifold, $\Omega\subseteq {\cal M}$ an open set and
  $X:{\Omega}\to T{\cal M}$   a continuously differentiable
 vector field. Let $p_*\in \Omega$, $R>0$  and $\kappa:=\sup\{t\in [0, R): B_{ t}(p_*)\subset \Omega\}$.
  Suppose that $X(p_*)=0$,  $\nabla X (p_*)$ is invertible and  there exists an $f:[0,\; R)\to \mathbb{R}$ continuously differentiable such that
   \begin{equation}\label{Hyp:MH}
\left\|\nabla X (p_*)^{-1}[P_{\zeta,1,0} \,\nabla X (p) -  \,P_{\zeta,\tau,0}\nabla X (\zeta(\tau))P_{\zeta,1,\tau}]\right\|
\leq f'\left(d(p_*, p) \right)-f'\left(\tau d(p_*, p)\right),
  \end{equation}
  for all $\tau \in [0,1]$, $p\in B_{\kappa}(p_*)$, where $\zeta:[0, 1]\to {\cal M}$ is a minimizing geodesic from $p_*$ to $p$ and
\begin{itemize}
  \item[{\bf h1)}]  $f(0)=0$ and $f'(0)=-1$;
  \item[{\bf  h2)}]  $f'$ is  strictly increasing.
\end{itemize}
  Let  $0\leq \vartheta<1/K_{p_*}$,  $\nu:=\sup\{t\in [0, R): f'(t)< 0\},$ $\rho:=\sup\{\delta\in (0, \nu): [(1+\vartheta)|t-f(t)/f'(t)|/t+\vartheta ]<1/K_{p_*},\, t\in (0, \delta)\} $  and
  $$
r:=\min \left\{\kappa, \,\rho,\, r_{p_*} \right\}.
  $$
Then the sequence  generated by the Inexact Newton method for solving $X(p)=0$ with starting point $p_0\in B_{r}(p_*)\setminus\{p_*\}$  and residual relative error tolerance $ \theta$,
\begin{equation} \label{eq:DNS}
p_{k+1} =\exp_{p_k}\left(S_k \right), \qquad \|  X(p_k) +  \nabla X ({p_k})S_k\| \leq \theta \|X ({p_k}) \|,  \qquad k=0,1,\ldots\,,
\end{equation}
\begin{equation} \label{eq:ltheta}
 0\leq  \mbox{cond}(\nabla X(p_*)) \theta \leq \vartheta/ \left[2/{|f'(d(p_*,p_0))|}-1 \right],
  \end{equation}
is well defined (for any particular choice of each $S_k\in T_{p_k}M$), the sequence $\{p_k\}$ is contained in $B_{r}(p_*)$ and  converges to the point $p_*$ which is the unique zero of $X$ in $B_{\sigma}(p_*)$, where $\sigma:=\sup\{t\in (0, \kappa): f(t)< 0\}$,  and we have that:
  \begin{equation}
    \label{eq:q2}
    d(p_*,\,p_{k+1}) \leq K_{p_*} \left[(1+\vartheta)  \frac{\left|  d({p_*},p_k) -\displaystyle \frac{f(d({p_*},p_k))}{f'(d({p_*},p_k))} \right|}{d({p_*},p_k)}+\vartheta \right] d(p_*,\,p_k), \qquad k=0,1,\ldots\,,
  \end{equation}
and  $\{p_k\}$ converges linearly to $p_*$. If, in additional, the function $f$ satisfies the  following condition
  \begin{itemize}
  \item[{\bf  h3)}]  $f'$ is  convex,
\end{itemize}
 then   there holds
 \begin{equation} \label{eq:q4}
    d(p_*,\,p_{k+1}) \leq K_{p_*} \left[(1+\vartheta)  \frac{\left|  d({p_*},p_0) -\displaystyle \frac{f(d({p_*},p_0))}{f'(d({p_*},p_0))} \right|}{d^2({p_*},p_0)} d({p_*},p_k)+\vartheta \right] d(p_*,\,p_k), \quad k=0,1,\ldots\,.
  \end{equation}
 as a consequence,  the sequence $\{p_k\}$ converges to  $p_*$     with linear rate  as follows
\begin{equation} \label{eq:q3}
    d(p_*,\,p_{k+1}) \leq K_{p_*} \left[(1+\vartheta)  \frac{\left|  d({p_*},p_0) -\displaystyle \frac{f(d({p_*},p_0))}{f'(d({p_*},p_0))} \right|}{d({p_*},p_0)}+\vartheta \right] d(p_*,\,p_k), \quad k=0,1,\ldots\,.
  \end{equation}
  \end{theorem}
  \begin{remark}
 First note that from  simple algebraic manipulation we have the following equality
 $$
 \frac{\left|  d({p_*},p_k) -\displaystyle \frac{f(d({p_*},p_k))}{f'(d({p_*},p_k))} \right|}{d({p_*},p_k)} = \left | 1- \frac{1}{f'(d({p_*},p_k))} \frac{f(d({p_*},p_k))-f(0)}{d({p_*},p_k)-0}\right|.
 $$
Since  the sequence $\{p_k\}$ is contained in $B_{r}(p_*)$ and  converges to the point $p_*$ then it is easy to see that right hand side   of last equality goes to zero  as $k$ goes to infinity.  Therefore in Theorem~\ref{th:nt} if taking $\vartheta=\vartheta_k$   in each iteration
and letting $\vartheta_k$ goes to zero (in this case,  $\theta=\theta_k$ also goes to zero)  as $k$ goes to infinity, then  \eqref{eq:q2}  implies that  $\{p_k\}$
converges to $p_*$  with asymptotic  superlinear rate.

Note that letting   $\vartheta =0$ in Theorem~\ref{th:nt} which implies from \eqref{eq:ltheta} that $\theta=0$, the linear  equation  in \eqref{eq:DNS} is solved exactly. Therefore  \eqref{eq:q4} implies that $\{p_k\}$ converges to $p_*$ with quadratic rate.
  \end{remark}
From now on, we assume that the hypotheses of Theorem~\ref{th:nt}  hold with the exception of {\bf h3}, which will be considered to hold only when explicitly  stated.
\section{Preliminary results} \label{sec:pr}
The scalar function $f$ in Theorem~\ref{th:nt}  is called a {\it majorant function} for vector field $X$ at a point $p_*$.    In this section we analyze some basic properties of $f$ and the main relationships between   $f$ and  $X$.
\subsection{The majorant function} \label{sec:PMF}
We begin by proving that the constants   $\kappa$, $\nu$ and $\sigma$ are positives.
\begin{proposition}  \label{pr:incr1}
The constants $ \kappa,\, \nu $ and $\sigma$ are positives and $t-f(t)/f'(t)<0$ for all $t\in (0,\,\nu).$
\end{proposition}
\begin{proof}
Since $\Omega$ is open and $p_*\in \Omega$, we  conclude that $\kappa>0$. As $f'$ is continuous in $0$ with  $f'(0)=-1$, there exists $\delta>0$ such that $f'(t)<0$ for all $t\in (0,\, \delta)$, so  $\nu>0$. Because $f(0)=0$ and $f'$ is continuous in $0$ with  $f'(0)=-1$,  there exists $\delta>0$ such that $f(t)<0$ for all $t\in (0, \delta)$, hence $\sigma>0$.

Assumption {\bf  h2} implies that $f$ is strictly convex,  so using the strict convexity of $f$ and the first equality in assumption {\bf  h1} we have
$
f(t)-tf'(t)<f(0)=0
$
for all  $t\in  (0,\, R).$
If $t\in (0, \,\nu)$ then $f'(t)<0$, which combined with the last inequality yields the desired inequality.
\end{proof}

According to {\bf h2} and  definition of $\nu$, we have  $f'(t)< 0$ for all
$t\in[0, \,\nu)$.  Therefore Newton iteration map for  $f$ is well defined in $[0,\, \nu)$. Let us call it $n_f$,
\begin{equation} \label{eq:def.nf}
  \begin{array}{rcl}
  n_f:[0,\, \nu)&\to& (-\infty, \, 0],\\
    t&\mapsto& t-f(t)/f'(t).
  \end{array}
\end{equation}
Because  $f'(t)\neq 0$ for all $t\in[0, \,\nu)$  the Newton iteration map $n_f$ is a continuous function.
\begin{proposition}  \label{pr:incr3}
$
\lim_{t\to 0}|n_f(t)|/t=0.
$
As a consequence  $\rho>0 $ and
$(1+\vartheta)|n_f(t)|/t+\vartheta <1/K_{p_*}$  for all $ t\in (0, \, \rho)$.
\end{proposition}
\begin{proof}
Using definition in  \eqref{eq:def.nf},  Proposition \ref{pr:incr1},  $f(0)=0$ and definition of $\nu$, a simple algebraic manipulation gives
\begin{equation} \label{eq:rho}
\frac{|n_f(t)|}{t}= \frac{f(t)/f'(t)-t}{t}=\frac{1}{f'(t)} \frac{f(t)-f(0)}{t-0}-1, \qquad t\in (0,\,\nu).
\end{equation}
Because  $f'(0)\neq 0$ the first statement follows by taking the limit in~\eqref{eq:rho} as $t$ goes to $0$.

Since $\lim_{t\to 0}|n_f(t)|/t=0$ and $\vartheta<1/K_{p_*}$ the first equality in  \eqref{eq:rho} implies that there exists $\delta>0$ such that
$$
(1+\vartheta)[f(t)/f'(t)-t ]/t+\vartheta <1/K_{p_*}, \qquad  \;  t\in (0, \delta).
$$
Therefore  from definition of $\rho$ and \eqref{eq:def.nf} the  last result of the  proposition follows.
\end{proof}
\begin{proposition}  \label{pr:sif}
 If  $f$ satisfies   {\bf  h3} then the function  $(0,\, \nu) \ni t \mapsto n_f(t)|/t^2$   is  increasing.
\end{proposition}
\begin{proof}
Using definition of $n_f$ in \eqref{eq:def.nf}, Proposition~\ref{pr:incr1} and {\bf h1} we obtain, after simples algebraic manipulation,  that
\begin{equation} \label{eq:deta}
 \frac{|n_{f}(t)|}{t^2}=
\frac{1}{|f'(t)|}\int_{0}^{1}\frac{f'(t)-f'(\tau t)}{t} \,d
\tau, \qquad \forall \;t\in (0,\, \nu).
\end{equation}
On the other hand as $f'$ is strictly  increasing  the map $ [0,\, \nu) \ni t \mapsto [f'(t)-f'(\tau t)]/t$  is  positive for all $\tau
\in (0, 1)$. From {\bf  h3}   $f'$ is convex, so we conclude that the last map  is
increasing. Hence the second term in the right hand side of \eqref{eq:deta} is positive and increasing. Assumption   {\bf  h2} and definition of $\nu$ imply that  the first term in the right hand side of \eqref{eq:deta} is also positive and strictly  increasing. Therefore we conclude that  the left hand side of \eqref{eq:deta} is increasing and the statement of the proposition follows.
\end{proof}

\subsection{Relationship between the majorant function and the vector field} \label{sec:MFNLO}
We  present the main  relationships between the majorant function $f$ and the vector field $X$.
\begin{lemma} \label{wdns}
Let  $p\in \Omega\subseteq {\cal M}$. If \,\,$d(p_*, p)<\min\{\kappa, \nu\}$  then $\nabla X(p)$ is invertible and
$$
\|\nabla X(p)^{-1}P_{\zeta, 0,1}\nabla X(p_*)\|\leq  1/|f'(d(p_*, p))|
$$
where $\zeta:[0, 1]\to {\cal M}$  is a minimizing geodesic from $p_*$ to $p$.
In particular $\nabla X(p)$ is invertible for all $p\in B_{ r}(p_*)$ where $r$ is as defined in Theorem \ref{th:nt}.
\end{lemma}
\begin{proof}
See Lemma 4.4 of \cite{FS2012}.
\end{proof}
\begin{lemma} \label{l:theta}
Let $p\in \Omega\subseteq {\cal M}$.   If \,\,$d(p_*, p)\leq d(p_*, p_0) <\min\{\kappa, \nu\}$,  then there holds
  $$
     \mbox{cond}(\nabla X(p)) \leq \mbox{cond}(\nabla X(p_*)) \left[ 2/|f'(d(p_*,p_0))|-1\right].
  $$
   As a consequence,   $\theta \mbox{cond}(\nabla X(p)) \leq \vartheta .$
\end{lemma}
\begin{proof}
Let $I_{p_*}:T_{p_*}{\cal M}\to T_{p_*}{\cal M}$ the identity operator, $p\in B_{ \kappa}(p_*)$ and $\zeta:[0, 1]\to {\cal M}$ a minimizing geodesic from $p_*$ to $p$.
Since $P_{\zeta, 0,0}=I_{p_*}$ and $P_{\zeta, 0,1}$ is an isometry  we obtain 
  \begin{align*}
    \left\|\nabla X(p_*)^{-1}P_{\zeta, 1,0}\nabla X(p)P_{\zeta, 0,1}-I_{p_*}\right\|= \left\|\nabla X(p_*)^{-1}[P_{\zeta, 1,0}\nabla X(p)-P_{\zeta, 0,0}\nabla X(p_*)P_{\zeta, 1,0}]\right\|.
  \end{align*}
As $d(p_*, p)<\nu$ we have $f'(d(p_*, p))<0$.  Using  the last equation,  \eqref{Hyp:MH} and {\bf h1} we conclude that
$$
\|\nabla X(p_*)^{-1}P_{\zeta, 1,0}\nabla X(p)P_{\zeta, 0,1}-I_{p_*}\|\leq f'(d(p_*, p))+1.
$$
Since $P_{\zeta, 0,1}$ is an isometry and  $\|\nabla X(p)\| \leq \|\nabla X(p_*)\|\|\nabla X(p_*)^{-1}P_{\zeta, 1,0}\nabla X(p)P_{\zeta, 0,1}\|$,  triangular   inequality together with  above inequality  imply
$$
\|\nabla X(p)\| \leq  \|\nabla X(p_*)\| \left[f'(d(p_*, p))+2\right].
$$
On the other hand, it is easy  to see from Lemma~\ref{wdns} that $\|\nabla X(p)^{-1}\| \leq \|\nabla X(p_*)^{-1}\|/|f'(d(p_*,p))| $. Therefore, combining two last inequalities and definition of condition number we obtain
 $$
     \mbox{cond}(\nabla X(p)) \leq \mbox{cond}(\nabla X(p_*)) \left[ 2/|f'(d(p_*,p))|-1\right].
  $$
Since $f'$ is strictly increasing, $f'<0$ in $[0, \nu)$ and $d(p_*, p)\leq d(p_*, p_0) <\min\{\kappa, \nu\}$,  the first inequality of the lemma follows from last inequality.

The last inequality of the lemma follows from  \eqref{eq:ltheta} and first inequality.
\end{proof}
 The linearization error of $X$  at a point in $B_{\kappa}(p_*)$ is defined by:
\begin{equation}\label{eq:def.er}
  E_X( p_*, p):= X(p_*)-P_{\alpha, 0,1}\left[ X(p)+\nabla X(p)\alpha'(0)\right],\qquad p\, \in B_{\kappa}(p_*),
\end{equation}
where $\alpha:[0, 1]\to {\cal M}$ is a minimizing geodesic from $p$ to $p*$.
We will bound this error by the error in the linearization on the
majorant function $f$,
\begin{equation}\label{eq:def.erf}
        e_f(t,u):= f(u)-\left[ f(t)+f'(t)(u-t)\right],\qquad t,\,u \in [0,R).
\end{equation}
\begin{lemma}\label{pr:taylor}
Let $p\in \Omega\subseteq {\cal M}$.  If $d(p_*, p)\leq \kappa$ then $\|\nabla X (p_*)^{-1}  E_X(p_*, p) \|\leq e_f(d(p_*, p),0)$.
\end{lemma}
\begin{proof}
See Lemma 4.5 of \cite{FS2012}.
\end{proof}

\begin{lemma} \label{le:bns}
Let $p\in \Omega\subseteq {\cal M}$.  If   $d(p_*, p)< r$  then
$$
\left\|\nabla X ({p}) ^{-1} X(p)\right \| \leq \frac{f(d(p_*,\,p))}{f'(d(p_*,\,p))}, \qquad p\in B_{r}(p_*).
$$
\end{lemma}
\begin{proof}
Since $X(p_*)=0$, the inequality is trivial for $p=p_*$. Now assume that  $0<d(p_*, p)< r$.
Lemma~\ref{wdns} implies that  $\nabla X(p) $ is invertible. Let $\alpha:[0, 1]\to {\cal M}$ be a minimizing geodesic from $p$ to $p_*$. Because $X(p_*)=0$, the definition of $E_X(p_*, p)$ in \eqref{eq:def.er} and direct  manipulation yields
$$
-\nabla X(p)^{-1}P_{\alpha, 1,0}E_X(p,p_*)= \nabla X(p)^{-1}X(p)+\alpha'(0).
$$
Using  the above equation,  Lemma \ref{wdns} and Lemma \ref{pr:taylor}, it is easy to conclude that
\begin{align*}
  \|\nabla X(p)^{-1}X(p)+\alpha'(0)\|&\leq\| -\nabla X(p)^{-1}P_{\alpha, 1,0}\nabla X(p_*)\|\|\nabla X(p_*)^{-1}E_F(p,p_*)\|\\
                                    &\leq e_f(d(p_*, p), 0)/|f'(d(p_*, p))|.
\end{align*}
As  $f(0)=0$,   definition of $e_f$  gives
$
e_f(d(p_*, p), 0)/|f'(d(p_*, p))|=-d(p_*, p) +f(d(p_*, p))/f'(d(p_*, p)),
$
which combined  with  last inequality  yields 
$$
 \|\nabla X(p)^{-1}X(p)+\alpha'(0)\|\leq -d(p_*, p) +f(d(p_*, p))/f'(d(p_*, p)).
$$
Since $\|\alpha'(0)\|=d(p_*, p)$, after simples algebraic manipulation we conclude
$$
 \|\nabla X(p)^{-1}X(p)\|\leq \|\nabla X(p)^{-1}X(p)+\alpha'(0)\| +d(p_*, p),
$$
which combined with  last   inequality yields  the desired result.
\end{proof}
The outcome of an Inexact Newton iteration is any point satisfying
some error tolerance. Hence, instead of a mapping for Newton
iteration, we shall deal with a \emph{family} of mappings describing
all possible inexact iterations.
\begin{definition} \label{def:inire}
  For $0\leq \theta$, $\mathcal{N}_\theta$ is the family of maps
  $N_\theta:B_{r}(p_*)\to  \banacha$ such that
  \begin{equation}
    \label{eq:in.map}
    \left\|X(p)+\nabla X(p)\exp_{p}^{-1}N_\theta(p)\right\|\leq \theta \left\|X(p)\right\|,  \qquad p\in B_{r}(p_*).
  \end{equation}
\end{definition}
If $p\in B_{r}(p_*)$ then $\nabla X(p)$ is
non-singular. Therefore for $\theta= 0$ the family $\mathcal{N}_0$ has a
single element,
namely,  the exact Newton iteration map
\begin{equation} \label{NF}
  \begin{array}{rcl}
  N_{0}:B_{r}(p_*) &\to& {\cal M}\\
    p&\mapsto&\exp_{p}\left(- \nabla X ({p}) ^{-1} X(p)\right).
  \end{array}
\end{equation}
Trivially, if $0\leq \theta \leq \theta '$ then
$\mathcal{N}_0\subset\mathcal{N}_\theta\subset\mathcal{N}_{\theta '}$.
Hence $\mathcal{N}_\theta$ is non-empty for all $\theta\geq 0$.
\begin{remark}
  For any $\theta \in (0,1)$ and $N_\theta\in\mathcal{N}_\theta$
  \[ N_\theta(p)=p\iff X(p)=0,\qquad p\in B_{r}(p_*).
  \]
  This means that the fixed points of the Inexact Newton iteration
  $N_\theta$ are the same fixed points of the \emph{exact} Newton
  iteration, namely,  the zeros of $X$.
\end{remark}
\begin{lemma} \label{le:cl}
Let  $\theta $ be such that $0\leq \theta  \mbox{cond}(\nabla X(p_*))  \leq \vartheta/ \left[1 + 2/{|f'(d(p_*,p_0))|} \right]$ and $p\in \Omega\subseteq {\cal M}$. If   $d(p_*, p)\leq d(p_*, p_0) <r$ and  $N_\theta \in \mathcal{N}_\theta$ then
$$
d(p_*,N_\theta(p))\leq   K_{p_*}\left[(1+\vartheta)\frac{|n_f(d({p_*},p))|}{d({p_*},p)}+\vartheta\right]\,d(p_*,\,p), \qquad p\in B_{r}(p_*).
$$
As a consequence,   $N_{\theta}(B_{r}(p_*))\subset B_{r}(p_*). $
\end{lemma}
\begin{proof}
Since $X(p_*)=0$, the inequality is trivial for $p=p_*$. Now, assume that  $0<d(p_*, p)\leq r$.  Let $\alpha:[0, 1]\to {\cal M}$ be a minimizing geodesic from $p$ to $p_*$. After  simple  algebraic manipulations, triangular inequality and definition of the linearization error we obtain
\begin{multline} \label{eq:cv1}
\| \exp_{p}^{-1}N_\theta(p)-\alpha'(0)\| \leq   \left\| \nabla X(p)^{-1}\left[\nabla X(p)\exp_{p}^{-1}N_\theta(p)+X(p)\right] \right\|+\left\|\nabla X(p)^{-1}E_X(p_*, p)\right \|.
\end{multline}
Using  Definition~\eqref {def:inire}  the first term in the right hand side of the above inequality is bounded by
$$
 \left\| \nabla X(p)^{-1}\left[\nabla X(p)\exp_{p}^{-1}N_\theta(p)+X(p)\right]\right\|\leq  \left\| \nabla X(p)^{-1} \right\| \theta \|X(p)\|.
$$
Now, since   $\|X(p)\| \leq \left\| \nabla X(p)\right\| \left\|\nabla X(p)^{-1}X(p)\right\|$ we obtain from Lemma~\eqref{le:bns} that
$$
|X(p)\| \leq\left\| \nabla X(p)\right\| \frac{f(d(p_*,\,p))}{f'(d(p_*,\,p))}.
$$
Definition of condition number and two above inequalities imply
\begin{equation} \label{eq:cv2}
 \left\| \nabla X(p)^{-1}\left[\nabla X(p)\exp_{p}^{-1}N_\theta(p)+X(p)\right]\right\|\leq \theta  \mbox{cond}(\nabla X(p)) \frac{f(d(p_*,\,p))}{f'(d(p_*,\,p))}.
\end{equation}
Now, combining  Lemma~\eqref{pr:taylor} and Lemma~\eqref{wdns} the second term in \eqref{eq:cv1} is bounded by
$$
\left\|\nabla X(p)^{-1}E_X(p_*, p),\right \| \leq \dfrac{1}{|f'(d(p_*,p))|}e_f(d(p_*,p),0).
$$
Therefore, \eqref{eq:cv1}, \eqref{eq:cv2} and last inequality give us
$$
\| \exp_{p}^{-1}N_\theta(p)-\alpha'(0)\| \leq \theta  \mbox{cond}(\nabla X(p)) \frac{f(d(p_*,\,p))}{f'(d(p_*,\,p))} +   \dfrac{1}{|f'(d(p_*,p))|}e_f(d(p_*,p),0).
$$
Since Lemma~\eqref{l:theta} implies $\theta  \mbox{cond}(\nabla X(p))\leq \vartheta$,  after simple algebraic manipulation and taking in account definitions of  $e_f$ and $n_f$ the  above inequaliy becomes 
$$
\| \exp_{p}^{-1}N_\theta(p)-\alpha'(0)\| \leq \left[(1+\vartheta)\frac{|n_f(d({p_*},p))|}{d({p_*},p)}+\vartheta\right]\,d(p_*,\,p).
$$
Note that, as $d(p_*, p)\leq r<\rho$,  second part of  Proposition \eqref{pr:incr3} implies that  the term in brackets  of last inequality is less than  $1/K_{p_*}\leq 1$. So left hand side of last inequality is less than $r\leq r_{p_*}$.  Therefore  letting  $p=p_*$,  $q=p$,  $v=\alpha'(0)$, $u = \exp_{p}^{-1}N_\theta(p)$  in   Definition \ref{def:kp}  we conclude that
$$
d(p_*,N_\theta(p))\leq K_{p_*}\| \exp_{p}^{-1}N_\theta(p)-\alpha'(0)\|.
$$
Finally  combining two above inequalities  the inequality of the lemma follows.

Take $p\in B_{r}(p_*)$.  Since $d(p_*, p)<r$ and $ r\leq \rho$, the first part of the lemma and the second part of Proposition~\ref{pr:incr3} imply that
$
d(p_*, N_X(p))<d(p_*, p)
$
and the  result  follows.
\end{proof}

\section{The Newton sequence} \label{sec:proof}
In this section  we prove  Theorem~\ref{th:nt}.  Let  $0\leq \theta$ satisfying \eqref{eq:ltheta} and  $N_{\theta}\in \mathcal{N}_\theta$, where $\mathcal{N}_\theta$  is defined in Definition~\ref{def:inire}.  Therefore   \eqref{eq:DNS} together with Definition~\ref{def:inire} implies that   the sequence $\{p_k\}$  satisfies
\begin{equation} \label{NFS}
p_{k+1}=N_{\theta}(p_k),\qquad k=0,1,\ldots \,,
\end{equation}
which is indeed an equivalent definition of this sequence.
\begin{proof}[\bf Proof of Theorem~\ref{th:nt}:]
Since $p_0\in B_{r}(p_*)$,  $r\leq \nu$ and $0<\theta  \mbox{cond}(\nabla X(p_*))  \leq \vartheta/ \left[2/{|f'(d(p_*,p_0))|}-1 \right]$, combining \eqref{NFS},  the inclusion  $N_{\theta}(B_{r}(p_*)) \subset B_{r}(p_*)$ in Lemma~\ref{le:cl} and Lemma~\ref{wdns}, it is easy to conclude that by an induction argument the sequence  $\{p_k\}$ is well defined and remains  in $B_{r}(p_*)$.

Now we are going to prove that $\{p_k \}$ converges towards $p_*$. Since $d(p_*,p_{k})<r$,  for $ k=0,1,\ldots \,$,  we obtain from  \eqref{NFS} and  Lemma~\ref{le:cl}  that
\begin{equation}\label{eq:conv1.0}
d(p_*,\,p_{k+1})\leq   K_{p_*}\left[(1+\vartheta)\frac{|n_f(d({p_*},p_k))|}{d({p_*},p_k)}+\vartheta\right]\,d(p_*,\,p_k).
\end{equation}
As  $d(p_*,p_{k})<r\leq \rho$,  for $ k=0,1,\ldots \,$,   using  second statement in Proposition~\ref{pr:incr3} and last inequality we conclude that $0\leq d(p_*,\,p_{k+1})<d(p_*,\,p_{k})$, for $ k=0,1,\ldots \,$.   So $\{d(p_*,p_{k})\}$ is strictly decreasing  and bounded below which implies that it   converges. Let  $\ell_*:=\lim_{k\to \infty}d(p_*,p_{k})$. Because  $\{d(p_*,p_{k})\}$ rests in $(0, \,\rho)$ and is strictly decreasing we have $0\leq \ell_*<\rho$. We are going to show that  $\ell_{*}=0$.  If $0<\ell_{*}$ then  letting $k$ goes to infinity in \eqref{eq:conv1.0},    the  continuity of $n_f$ in $[0, \rho)$   and  Proposition~\ref{pr:incr3} imply that
\begin{equation}\label{eq:conv1.1}
\ell_{*} \leq   K_{p_*}\left[(1+\vartheta)\frac{|n_f(\ell_{*})|}{\ell_{*}}+\vartheta\right]\,\ell_{*}<\ell_{*},
\end{equation}
which is an absurd. Hence we must have $\ell_{*}=0$.  Therefore  the convergence of $\{p_k \}$ to $p_*$ is proved.   The uniqueness of $p_*$ in $B_{ \sigma}(p_*)$ was proved in Lemma 5.1 of \cite{FS2012}.

For  proving   the  equality in  \eqref{eq:q2}  it is sufficient to use equation \eqref{eq:conv1.0} and definition of $n_f$ in \eqref{eq:def.nf}.   As $d(p_*,p_{k})<r\leq \rho$,  for $ k=0,1,\ldots \,$,  $\lim_{k\to \infty} d(p_*,p_{k})=0$ and  by hypothesis  $\vartheta<1/K_{p_*}$ thus using definition of $n_f$ and first statement in  Proposition~\ref{pr:incr3} we conclude
$$
\lim_{k\to \infty}K_{p_*} \left[(1+\vartheta)\frac{\left| d({p_*},p_k)-\displaystyle \frac{f(d({p_*},p_k))}{f'(d({p_*},p_k))}\right|}{d({p_*},p_k))} + \vartheta\right]= K_{p_*}\vartheta<1.
$$
 which implies  the linear convergence of   $\{p_k \}$ to $p_*$  in \eqref{eq:q2}.

 Now we are going to prove the inequality in \eqref{eq:q4}:  If  $f$ satisfies   {\bf  h3}  then using  definition of $n_f$ and Proposition~\ref{pr:sif} we conclude
 $$
(1+\vartheta)\frac{\left| d({p_*},p_k)-\displaystyle \frac{f(d({p_*},p_k))}{f'(d({p_*},p_k))}\right|}{d^2({p_*},p_k))} d({p_*},p_k)) + \vartheta \leq (1+\vartheta)\frac{\left| d({p_*},p_0)-\displaystyle \frac{f(d({p_*},p_0))}{f'(d({p_*},p_0))}\right|}{d^2({p_*},p_0))} d({p_*},p_k))+ \vartheta.
$$
  As the quantity of the left hand side of the last inequality is equal to quantity in the brackets  of   \eqref{eq:q2},   the inequality in  \eqref{eq:q4}   follows from  \eqref{eq:q2} and last inequality.
  
Since        $\{d(p_*,p_{k})\}$ is strictly decreasing,    the inequality in \eqref{eq:q3} follows from \eqref{eq:q4} and we  conclude the proof of the theorem.
\end{proof}
\section{Special Cases} \label{apl}
In this section, we present two   special cases of Theorem~\ref{th:nt}.
\subsection{Convergence result  under H\"{o}lder-like  condition}
For null error tolerance, the next theorem on  Inexact Newton's method under a H\"{o}lder-like condition merges in  Theorem~ 7.1  of \cite{FS2012}. 
\begin{theorem} \label{th:HV}
Let $\cal M$ be a Riemannian manifold, $\Omega\subseteq {\cal M}$ an  open set  and
  $X:{\Omega}\to T{\cal M}$ a  continuously differentiable
 vector field.
 Take $p_*\in \Omega$, $R>0$  and let $\kappa:=\sup\{t\in [0, R): B_{ t}(p_*)\subset \Omega\}$.
  Suppose that $X(p_*)=0$,  $\nabla X (p_*)$ is invertible and there exist  constants  $L>0$ and $0\leq \mu< 1$   such that
  \begin{equation}\label{Hyp:HC}
\left\|\nabla X (p_*)^{-1}[P_{\zeta,1,0} \,\nabla X (p) -  \,P_{\zeta,\tau,0}\nabla X (\zeta(\tau))P_{\zeta,1,\tau}]\right\|
\leq L(1-\tau^{\mu})d(p_{*}, p)^{\mu},
  \end{equation}
  for all $\tau \in [0,1]$ and  $p\in B_{\kappa}(p_*)$, where $\zeta:[0, 1]\to {\cal M}$ is a minimizing geodesic from $p_{*}$ to $p$.
Let $r_{p_*}$  be the injectivity radius of $\cal{M}$ in $p_*$, $K_{p_*}$ as  in Definition \ref{def:kp}, $0 \leq \vartheta < 1/K_{p_*}$ and
$$
r:=\min \left\{\kappa, \,\left[(\mu+1)\big/\left(L\left(\frac{1+K_{p_*}}{1-K_{p_*}\vartheta}\mu+1\right)\right)\right]^{1/\mu}, \, r_{p_*} \right\}.
$$
Then the sequence generated by the Inexact Newton method for solving $X(p)=0$ with starting point $p_0\in B_{ r}(p_*)\setminus\{p_*\}$ and residual relative error tolerance $\theta$,
\begin{equation} \label{eq:DNS2}
p_{k+1} =\exp_{p_k}\left(S_k \right), \qquad \|  X(p_k) +  \nabla X ({p_k})S_k\| \leq \theta \|X ({p_k}) \|,  \qquad k=0,1,\ldots\,,
\end{equation}
\begin{equation} \label{eq:ltheta2}
 0\leq  \mbox{cond}(\nabla X(p_*)) \theta \leq \vartheta \; \frac{1+Ld(p_*, p_0)^{\mu}}{1-Ld(p_*, p_0)^{\mu}},
  \end{equation}
is well defined (for any particular choice of each $S_k\in T_{p_k}M$), the sequence $\{p_k\}$ is contained in $B_{r}(p_*)$ and  converges to the point $p_*$ which is the unique zero of $X$ in  $B_{[(\mu+1)/L]^{1/\mu}}(p_*)$ and we have that:
$$
 d(p_*,p_{k+1}) \leq K_{p_*}\left[   (1+\vartheta) \frac{\mu L d({p_*},p_k)^\mu}{(\mu + 1)\left[1-Ld(p_*,p_k)^\mu \right]} +\vartheta \right] d({p_*},p_k),   \qquad k=0,1,\ldots\,,
$$
and  $\{p_k\}$ converges linearly to $p_*$. If, in additional,  $\mu=1$    then   there holds
\begin{equation} \label{eq:qch}
 d(p_*,p_{k+1}) \leq K_{p_*}\left[   (1+\vartheta) \frac{L}{2\left[1-Ld(p_*, p_0)\right]}d({p_*},p_k)  +\vartheta \right] d({p_*},p_k)  \qquad k=0,1,\ldots\,
\end{equation} 
as a consequence,  the sequence $\{p_k\}$ converges to  $p_*$     with linear rate  as follows
$$
 d(p_*,p_{k+1}) \leq K_{p_*}\left[   (1+\vartheta) \frac{ L d(p_*,p_0)}{2\left[1-Ld(p_*, p_0) \right]} +\vartheta \right] d({p_*},p_k)  \qquad k=0,1,\ldots\,
$$
\end{theorem}
\begin{proof}
 We can  prove that  $X$, $p_*$ and $f:[0, +\infty)\to \mathbb{R}$, defined by
$
f(t)=Lt^{\mu+1}/(\mu+1)-t,
$
satisfy the inequality \eqref{Hyp:MH} and the conditions  {\bf h1} and {\bf h2}  in Theorem \ref{th:nt}.  Moreover,  if $\mu=1$ then $f$ satisfies condition  {\bf h3}. It is easy to see that  $\rho$,  $\nu$ and $\sigma$, as defined in Theorem \ref{th:nt}, satisfy
$$
\rho=\left[\frac{(\mu+1)}{L \displaystyle \left(\frac{1+K_{p_*}}{1-K_{p_*}\vartheta}\mu+1\right)}\right]^{1/\mu} \leq ~ \nu=\frac{1}{L ^{1/\mu}}\, , \qquad  \sigma={[(\mu+1)/L]^{1/\mu}}.
$$
Therefore, the result follows  by invoking Theorem~\ref{th:nt}.
\end{proof}
\begin{remark}
Note that if vector field $X$ is  Lipschitz with constant $L$ then it satisfies the condition  \eqref{Hyp:HC} with $\mu=1$.

We remark that letting   $\vartheta =0$ in Theorem~\ref{th:HV} which implies from \eqref{eq:ltheta2} that $\theta=0$, the linear  equation  in \eqref{eq:DNS2} is solved exactly. Therefore  \eqref{eq:qch} implies that if $\mu=1$ then  $\{p_k\}$ converges to $p_*$ with quadratic rate.
\end{remark}

\subsection{Convergence result under Smale's condition }
For null error tolerance, the next theorem on  Inexact Newton's method under  Smale's condition merges in Theorem~ 7.2  of  \cite{FS2012}.   We note that Theorem~ 7.2  of  \cite{FS2012}  extends to the Riemannian context  Theorem~1.1 of \cite{D2003} (see also  Theorem~3.1 of  \cite{Wang2006}) which generalizes to the Riemannian context Corollary of Proposition~3 on p.~195 of \cite{Smale1986}, see also Proposition 1 p.~157 and Remark 1 p.~158 of \cite{Blum1998}.  

\begin{theorem}\label{theo:Smale}
Let $\cal M$ be an analytic Riemannian manifold, $\Omega\subseteq {\cal M}$ an open set  and
  $X:{\Omega}\to T{\cal M}$ an analytic vector field. Take $p_*\in \Omega$, $R>0$  and let $\kappa:=\sup\{t\in [0, R): B_{t}(p_*)\subset \Omega\}$.
  Suppose that $X(p_*)=0$,  $\nabla X (p_*)$ is invertible and
\begin{equation} \label{eq:SmaleCond}
 \gamma := \sup _{ n > 1 }\left\| \frac
{\nabla X(p_*)^{-1}\nabla^n X(p_*)}{n !}\right\|^{1/(n-1)}<+\infty.
\end{equation}
Let $r_{p_*}$  be the injectivity radius of $\cal{M}$ in $p_*$, $K_{p_*}$ as  in Definition \ref{def:kp}, $0 \leq \vartheta < 1/K_{p_*}$ and
$$
r:=\min \left\{\kappa, \,\frac{K_{p*}(1-3\vartheta)+4-\sqrt{K_{p*}^2(1-6\vartheta+\vartheta^2)+8K_{p*}(1-\vartheta)+8}}{4\gamma(1-K_{p_*}\vartheta)}, \, r_{p_*} \right\}.
$$
Then the sequence generated by the Inexact Newton method for solving $X(p)=0$ with starting point $p_0\in B_{ r}(p_*)\setminus\{p_*\}$ and residual relative error tolerance $\theta$,
\begin{equation} \label{eq:DNS3}
p_{k+1} =\exp_{p_k}\left(S_k \right), \qquad \|  X(p_k) +  \nabla X ({p_k})S_k\| \leq \theta \|X ({p_k}) \|,  \qquad k=0,1,\ldots\,,
\end{equation}
\begin{equation} \label{eq:ltheta3}
 0\leq  \mbox{cond}(\nabla X(p_*)) \theta \leq \vartheta \left[2[1-\gamma d(p_*,p_0)]^2-1\right], 
  \end{equation}
is well defined (for any particular choice of each $S_k\in T_{p_k}M$), the sequence $\{p_k\}$ is contained in $B_{r}(p_*)$ and  converges to the point $p_*$ which is the unique zero of $X$ in  $B_{1/(2\gamma)}(p_*)$ and we have that:

\begin{equation} \label{eq:q4S1}
    d(p_*,\,p_{k+1}) \leq K_{p_*} \left[(1+\vartheta)  \frac{\gamma}{2\left[1-\gamma d({p_*},p_0)\right]^2-1} d({p_*},p_k)+\vartheta \right] d(p_*,\,p_k), \quad k=0,1,\ldots\,.
  \end{equation}
 as a consequence,  the sequence $\{p_k\}$ converges to  $p_*$     with linear rate  as follows
\begin{equation} \label{eq:q4S2}
    d(p_*,\,p_{k+1}) \leq K_{p_*} \left[(1+\vartheta)  \frac{\gamma d({p_*},p_0)}{2\left[1-\gamma d({p_*},p_0)\right]^2-1}+\vartheta \right] d(p_*,\,p_k), \quad k=0,1,\ldots\,.
  \end{equation}

 \end{theorem}
We need the following results to prove the above theorem.
\begin{lemma} \label{lemma:qc1}
Let $\cal M$ be an analytic Riemannian manifold, $\Omega\subseteq {\cal M}$ an open set and
  $X:{\Omega}\to T{\cal M}$ an analytic vector field.   Suppose that
$p_*\in \Omega$,  $\nabla X (p_*)$ is invertible, $\gamma<+\infty$ and that $B_{
1/\gamma}(p_{*}) \subset \Omega$, where $\gamma$ is defined in
\eqref{eq:SmaleCond}. Then, for all $p\in B_{1/\gamma}(p_{*}),$
$$
\|\nabla X (p_*)^{-1}P_{\zeta, 1,0}\nabla^2 X(p))\| \leq  (2\gamma)/( 1- \gamma
d(p_*,p))^3,
$$
where $\zeta:[0, 1]\to {\cal M}$  is a minimizing geodesic from $p_*$ to $p$.
\end{lemma}
\begin{proof}
The proof follows the pattern  of Lemma~5.3 of \cite{Alvarez2008}.
\end{proof}
The next result is the Lemma 7.4 of \cite{FS2012}, it  gives an alternative condition  for  checking  condition \eqref{Hyp:MH}, whenever the vector field under consideration is twice continuously differentiable.
\begin{lemma} \label{lc}
Let $\cal M$ be an analytic Riemannian manifold, $\Omega\subseteq {\cal M}$ an open set and
  $X:{\Omega}\to T{\cal M}$ an analytic vector field.   Suppose that
$p_*\in \Omega$ and   $\nabla X (p_*)$ is invertible.  If there exists an \mbox{$f:[0,R)\to \mathbb {R}$} twice continuously differentiable such that
 \begin{equation} \label{eq:lc2}
\|\nabla X (p_*)^{-1}P_{\alpha, 1,0}\nabla^2 X(q))\|\leqslant f''(d(p_*,q)),\qquad \forall\,  q\in B_{\kappa}(p_*),
\end{equation}
where $\alpha:[0, 1]\to {\cal M}$  is a minimizing geodesic from $p_*$ to $q$, then $X$ and $f$
satisfy \eqref{Hyp:MH}.
\end{lemma}
\begin{corollary} \label{eq:ast}
Let $\cal M$ be an analytic Riemannian manifold, $\Omega\subseteq {\cal M}$ an open set and
  $X:{\Omega}\to T{\cal M}$ an analytic vector field. Take $p_*\in \Omega$  and let $\kappa:=\sup\{t\in [0, R): B_{t}(p_*)\subset \Omega\}$ and $\gamma<+\infty$ be as defined in \eqref{eq:SmaleCond}. Suppose that   $\nabla X (p_*)$ is invertible. Then
\begin{align*}
 \left\|\nabla X (p_*)^{-1}[P_{\zeta,1,0} \,\nabla X (p) -  \,P_{\zeta,\tau,0}\nabla X (\zeta(\tau))P_{\zeta,1,\tau}]\right\|
&\leq \frac{1}{( 1- \gamma d(p_*,p))^2}-\frac{1}{( 1- \tau\gamma d(p_*,p))^2}
\end{align*}
for all $\tau \in [0,1]$, $p\in B_{1/\gamma}(p_*)$, where $\zeta:[0, 1]\to {\cal M}$ a minimizing geodesic from $p_*$ to $p$
\end{corollary}
\begin{proof}
The proof follows by a combination  of  Lemma~\ref{lc} with  Lemma~\ref{lemma:qc1}.
\end{proof}

\noindent
{\bf [Proof of Theorem \ref{theo:Smale}]}. 
Assume that all hypotheses of   Theorem \ref{theo:Smale} hold. Consider the
real analytical  function $f:[0,1/\gamma) \to \mathbb{R}$ defined by
$$
f(t)=\frac{t}{1-\gamma t}-2t.
$$
It is straightforward to show that $f$ is  analytic and that
$$
f(0)=0, \quad f'(t)=1/(1-\gamma t)^2-2, \quad f'(0)=-1, \quad
f''(t)=(2\gamma)/(1-\gamma t)^3, \quad f'''(t)=6\gamma^2/(1-\gamma t)^4.
$$
It follows from  the last equalities that
$f$ satisfies {\bf h1},  {\bf h2} and {\bf h3}. Now, since $f'(t)=1/(1-\gamma t)^2-2$  we
conclude  from Corollary~\ref{eq:ast} that  $X$  and $f$ satisfy  \eqref{Hyp:MH} with $R=1/\gamma$. In this case, it is easy to see that the constants $\nu$,  $\rho$ and  $\sigma$, as defined in Theorem~\ref{th:nt},  satisfy
$$
\rho=\frac{K_{p*}(1-3\vartheta)+4-\sqrt{K_{p*}^2(\vartheta^2-6\vartheta+1)+8K_{p*}(1-\vartheta)+8}}{4\gamma(1-K_{p_*}\vartheta)}\leq \nu=\frac{\sqrt{2}-1}{\gamma\sqrt{2}},
$$
$\sigma=1/(2\gamma)$  and $f(0)=f(1/(2\gamma))=0$ and $f(t)<0$ for all $t\in (0,\, 1/(2\gamma))$.
Therefore, the result follows  by invoking Theorem~\ref{th:nt}.
\qed

\begin{remark}
We remark that letting   $\vartheta =0$ in Theorem~\ref{theo:Smale} which implies from \eqref{eq:ltheta3} that $\theta=0$, the linear  equation  in \eqref{eq:DNS3} is solved exactly. Therefore  \eqref{eq:q4S1} implies that  $\{p_k\}$ converges to $p_*$ with quadratic rate.
\end{remark}
\section{Final remarks } \label{rf}
The results in Theorem~\ref{th:nt} are dependent on the injective radius of the exponential map. It would  be  interesting to establish the convergence radius independent of the  injective radius of the exponential map.



\end{document}